\documentclass{article}[12pt]

\usepackage{amsmath}
\usepackage{amssymb}
\usepackage{amsthm}
\usepackage{enumerate}
\usepackage[pdftex]{hyperref} 

\newtheorem{theorem}{Theorem}
\newtheorem{lemma}[theorem]{Lemma}
\newtheorem{cor}[theorem]{Corollary}

\title{\bf
On a Diagonal Conjecture\\for
Classical Ramsey Numbers}

\author{Meilian Liang\\[-0.1ex]
\small School of Mathematics and Information Science\\[-0.6ex]
\small Guangxi University, Nanning 530004, P.R. China\\[-0.6ex]
\small {\tt gxulml@163.com}\\[1.3ex]\and
Stanis{\l}aw Radziszowski\\[-0.1ex]
\small Department of Computer Science\\[-0.6ex]
\small Rochester Institute of Technology, Rochester, NY 14623\\[-0.6ex]
\small {\tt spr@cs.rit.edu}\\[1.3ex]\and
Xiaodong Xu\\[-0.1ex]
\small Guangxi Academy of Sciences\\[-0.6ex]
\small Nanning 530007, P.R. China\\[-0.6ex]
\small {\tt xxdmaths@sina.com}\\[3.3ex]
}

\date{\today}

\begin{document}
\maketitle
\thispagestyle{empty}

\begin{abstract}
Let $R(k_1, \cdots, k_r)$ denote the classical $r$-color
Ramsey number for integers $k_i \ge 2$.
The Diagonal Conjecture (DC) for classical Ramsey numbers
poses that if $k_1, \cdots, k_r$ are integers no
smaller than 3 and $k_{r-1} \leq k_r$, then
$R(k_1, \cdots, k_{r-2}, k_{r-1}-1, k_r +1) \leq R(k_1, \cdots, k_r)$.
We obtain some implications of this conjecture, present
evidence for its validity, and discuss related problems.

Let $R_r(k)$ stand for the $r$-color Ramsey number $R(k, \cdots, k)$.
It is known that $\lim_{r \rightarrow \infty} R_r(3)^{1/r}$
exists, either finite or infinite, the latter
conjectured by Erd\H{o}s.
This limit is related to the Shannon capacity of
complements of $K_3$-free graphs.
We prove that if DC holds, and
$\lim_{r \rightarrow \infty} R_r(3)^{1/r}$
is finite, then
$\lim_{r \rightarrow \infty} R_r(k)^{1/r}$
is finite for every integer $k \geq 3$.
\end{abstract}

{\bf AMS classification subjects:} 05C55, 05C35

\bigskip
\section{Introduction}

Denote by $K_n$ the complete graph on
$n$ vertices. The classical multicolor Ramsey number
$R(k_1, \cdots, k_r)$ is the smallest positive integer $n$
such that if we color the edges of $K_n$ with $r$ colors,
then in this coloring there must be a monochromatic
$K_{k_i}$ whose all edges are in color $i$, for some
$i \in \{1, \dots, r\}$. In the diagonal case
$k=k_1 =  \cdots = k_r$ we will use the simpler
notation $R_r(k)=R(k_1, \cdots, k_r)$.

Wang Rui \cite{Wang} in a 2008 paper claimed to prove
that in the two-color case it holds that
$R(s,t)>R(s-1,t+1)$ for $s \leq t$.
Or, equivalently, one of his theorems states that as we
move away from the diagonal of the table with Ramsey numbers
$R(s,t)$, while preserving $s+t$, the values decrease.
Known values and bounds for Ramsey numbers \cite{Radziszowski}
do not contradict
this claim, and actually, it seems very plausible to be true.
Unfortunately, it is rather evident that its proof
in \cite{Wang} is not correct. The problems with this
paper are numerous, starting with a strange alternate
definition of Ramsey numbers, followed by unfounded
circular arguments between the alternate definitions.
Wang in his paper is addressing almost exclusively
two-color cases, but towards the end
he also makes some claims for more colors, though
again without what can be considered rigorous proofs.

We summarize the above in the following conjecture for
general multicolor Ramsey numbers, where two colors
are a special case.

\bigskip
\noindent
{\bf Diagonal Conjecture (DC).}\\
{\it
If $k_1, \cdots, k_r$
are integers no smaller than $3$, $r \ge 2$,
and $k_{r-1} \leq k_r$, then
$$R(k_1, \cdots, k_{r-2}, k_{r-1}-1, k_r +1) \leq R(k_1, \cdots, k_r).$$
}

If DC holds, then for the last two colors (and thus also
for any two fixed colors) as we move away from the diagonal,
while preserving $k_{r-1}+k_r$, the corresponding
Ramsey number cannot increase. We believe that a stronger
version of DC with $<$ instead of $\leq$ also holds.
Still, even the weaker version can be very hard to prove.

\medskip
In 1983,
Chung and Grinstead \cite{CG1983}
showed that
$\lim_{r \rightarrow \infty} R_r(3)^{1/r}$ exists,
though it is not known whether this limit is finite or infinite.
The same argument can be used to show that
$\lim_{r \rightarrow \infty} R_r(k)^{1/r}$ also exists
for all $k>3$, again finite or infinite.
Erd\H{o}s was inclined to think that
$\lim_{r \rightarrow \infty} R_r(3)^{1/r} = \infty$
(cf. \cite{Li,XRsur}).
This limit is also closely related to the Shannon capacity
of complements of $K_3$-free graphs (i.e. graphs with independence
number equal to 2), which was discussed
in an earlier paper by the second and third authors \cite{XRShannon}.

\medskip
Let $L_k=\lim_{r \rightarrow \infty} R_r(k)^{1/r}$.
By monotonicity of Ramsey numbers, we can easily
see that $L_{k+1} \ge L_k$ for all $k \ge 3$,
including the propagation of infinity to larger indices.
In this paper we obtain some consequences of the assumption
that the DC holds, we present evidence for its validity,
and discuss related problems. In particular, we prove
that if DC holds and
$\lim_{r \rightarrow \infty} R_r(3)^{1/r}$
is finite, then
$\lim_{r \rightarrow \infty} R_r(k)^{1/r}$
is finite for any integer $k \geq 3$. We also discuss
other relationships between DC and the sequence of $L_k$'s.

\bigskip
\section{Some Consequences of DC}

\begin{lemma}
If DC holds, then for every integer $k \geq 3$ we have
$$R_ {2r}(k) -1 \geq (R_r(k-1)-1)(R_r(k+1)-1).$$
\end{lemma}

\begin{proof}
An old result
obtained by Abbott \cite{Abb65}, also presented in \cite{XXER}
(Theorem 2, page 7), states that
if $k_j \ge 2$ for $1 \le j \le r$,
then for all $1 < i < r$ we have
$$R(k_1,\cdots,k_r) >
(R(k_1,\cdots,k_i)-1)(R(k_{i+1},\cdots,k_r)-1). \eqno{(1)}$$

\noindent
If DC holds, then we can apply it $r$ times
to $R_ {2r}(k)$ to obtain
$$R_ {2r}(k) -1 \geq R(k-1, \cdots,
k-1, k+1, \cdots, k+1) - 1.$$
Now, we can complete the proof
using inequality (1).
\end{proof}

\bigskip
\begin{theorem} \label{finite}
If DC holds and
$\lim_{r \rightarrow \infty} R_r(3)^{1 \over r}$
is finite, then
$\lim_{r \rightarrow \infty} R_r(k)^{1 \over r}$
is finite too,
for every integer $k \ge 3$.
\end{theorem}

\begin{proof}
For every integer $k \geq 3$, using Lemma 1 with DC, we have
$$(R_ {2r}(k) -1)^{1/r} \geq
(R_r(k-1)-1)^{1/r} (R_r(k+1)-1)^{1/r},$$
and thus

$${(R_ {2r}(k) -1)^{1\over{2r}} \over (R_r(k-1)-1)^{1/r} } \geq
{(R_r(k+1)-1)^{1/r} \over (R_ {2r}(k)-1)^{1\over{2r}}}.
\eqno{(2)}$$

\medskip
Clearly,
$\lim_{r \rightarrow \infty} R_r(i)^{1/r} =
\lim_{r \rightarrow \infty} (R_r(i)-1)^{1/r} $
for $i \in \{k-1, k, k+1\}$.
Taking it into account in inequality (2) leads to

$$\lim_{r \rightarrow \infty} {
{R_r(k)^{1 \over r} } \over {R_r(k-1)^{1 \over r}}  } \geq
\lim_{r \rightarrow \infty}{ {R_r(k+1)^{1 \over r} } \over
{R_r(k)^{1 \over r}}}.
\eqno{(3)}$$

\bigskip
Note that $R_r(2)=2$.
Finally, we can prove the claim of the theorem by induction on $k$.
The base case is for $k=3$, which is the finiteness of
$\lim_{r \rightarrow \infty} R_r(3)^{1 \over r}$.
The inductive step follows from the inequality (3).
\end{proof}

\bigskip
By Theorem \ref{finite}, we can see that if DC holds,
then either
$\lim_{r \rightarrow \infty} R_r(3)^{1 \over r} = \infty$,
or $\lim_{r \rightarrow \infty} R_r(k)^{1 \over r}$
is finite  for every $k \geq 3$. Or, equivalently,
DC and $\lim_{r \rightarrow \infty} R_r(k)^{1 \over r}  = \infty$
for any $k \ge 3$ implies that
$\lim_{r \rightarrow \infty} R_r(3)^{1 \over r}  = \infty$.
On the other hand, if
$\lim_{r \rightarrow \infty} R_r(3)^{1 \over r}$
were finite, then it would support our intuition
that the best known lower bounds for $R_r(3)$
are much closer to the exact values than the currently
best known upper bounds.

Table 1
presents the best known lower and upper bounds
on $R_r(3)$ for $r \leq 10$. The exact values
for $r=2,3$ are known, and it was conjectured that
$R_4(3)=51$, i.e. that the current lower bound
for $r=4$ is equal to the exact value \cite{XRsur}.
Lower bounds for higher $r$ in Table 1 are implied by
sum-free set constructions and related Schur
numbers (cf. \cite{XXER,XRsur}), in particular they
imply that
$\lim_{r \rightarrow \infty} R_r(3)^{1 \over r}
\geq 1073^{1 \over 6} \approx 3.1996$.

For the upper bound, a simple reasoning yields
$R_r(3) \leq 3r!$, while the best known general upper is just
a little better, namely, the third author et al.
proved that for $r \ge 4$ we have the bound
$R_r(3) \leq (e - {1 \over 6})r!+1 \approx 2.55r!$ \cite{rn3ub2002}.
The latter was proved based on the bound
$R_4(3) \leq 62$, which in turn was obtained with
the help of significant computations. This is the only case where
we know of an upper bound for a Ramsey number of this form
that is better than one obtained by simple steps using
smaller cases. Complete references
to lower and upper bounds and other general results
on $R_r(k)$ can be found in the dynamic survey paper
by the second author \cite{Radziszowski}.

If our perspective above that the lower
bounds in Table 1 are much closer to $R_r(k)$ than
the upper bounds is correct, it would add weight to
the case of $\lim_{r \rightarrow \infty} R_r(3)^{1 \over r}$
being finite, and thus by DC
and Theorem 2 also that all limits
$\lim_{r \rightarrow \infty} R_r(k)^{1 \over r}$ are finite.

\medskip
\begin{center}
\begin{tabular}{c|c|c}
$r$&lower bound&upper bound\cr
\hline
2&{\bf 6}&6\cr
3&{\bf 17}&17\cr
4&{\bf 51}&62\cr
5&162&307\cr
6&538&1838\cr
7&1682&12861\cr
8&5204&102882\cr
9&16146&925931\cr
10&51202&9259302\cr
\hline
\end{tabular}

\medskip
{\bf Table 1.} Known bounds on $R_r(3)$ for $r \le 10$.
\end{center}

\bigskip
We can prove the following Theorem \ref{largerthanone}
about the growth of the limits
$\lim_{r \rightarrow \infty} R_r(k)^{1 \over r}$
with increasing $k$, only assuming that DC holds.
However, we feel strongly that it also holds unconditionally.

\medskip
\begin{theorem} \label{largerthanone}
If DC holds, then for every integer $k \geq 3$, we have
$$\lim_{r \rightarrow \infty} {R_r(k)^{1 \over r}
\over {R_r(k-1)^{1 \over r}} }  > 1.$$
\end{theorem}

\medskip
\begin{proof}
Consider a general constructive lower bound for multicolor
Ramsey numbers $R_r(st+1) > (R_r(s+1)-1)(R_r(t+1)-1)$,
which can be obtained from a standard graph product construction
as described in \cite{XXER} (inequality (5) on page 4 there).
Using $s=k-1$ and $t=2$, it gives
$R_r(2k-1) > (R_r(3)-1)(R_r(k)-1)$.
We know that asymptotically $R_r(3)$ grows at least as
fast as $3.19^r$, but one can also easily observe
that $R_r(3) > 2^r$ holds for all $r$.
Thus
$${{R_r(2k-1) -1} \over {R_r(k)-1}} \geq R_r(3)-1 \ge 2^r,$$
and hence
$$ \lim_{r \rightarrow \infty} {R_r(2k)^{1 \over r} \over
{R_r(k)^{1 \over r}} }  \geq 2.\eqno{(4)}$$

Assume that DC holds, and for contradiction suppose
that for some $a \geq 3$ we have
$\lim_{r \rightarrow \infty} {R_r(a)^{1 \over r}/
{R_r(a-1)^{1 \over r}}}=1.$
Note that the inequality (3) in the proof of Theorem \ref{finite}
is a consequence of just DC, and it is valid for any $a \ge 3$.
Consequently, using (3) we can conclude that
$$\lim_{r \rightarrow \infty} {R_r(k)^{1 \over r}  \over
{R_r(k-1)^{1 \over r}}}=1$$
for any integer $k \geq a$.
This, however, leads to
$$\lim_{r \rightarrow \infty} {R_r(2a)^{1 \over r}  \over
{ R_r(a)^{1 \over r}} } =
\lim_{r \rightarrow \infty} \prod_{k=a+1} ^{2a}
{R_r(k)^{1 \over r}  \over { R_r(k-1)^{1 \over r}} } = 1,$$
which contradicts (4). This completes the proof of the theorem.
\end{proof}

\bigskip
\begin{cor}
For $k \ge 3$, let
$L_k=\lim_{r \rightarrow \infty} R_r(k)^{1/r}$,
and assume that DC holds. Then it is true that:

\smallskip
$(a)$ all $L_k$'s are finite or all of them are infinite, and

$(b)$ if $L_3$ is finite then $L_k < L_{k+1}$ for all $k \ge 3$.
\end{cor}

\smallskip
\begin{proof}
As discussed in the Introduction, all the limits $L_k$
exist and they satisfy $L_k \le L_{k+1}$, regardless of
whether DC holds or not. Thus, the claim $(a)$ follows from
Theorem \ref{finite} and claim $(b)$ follows from Theorem 
\ref{largerthanone}.
\end{proof}

We wish to note that clearly
$\lim_{k \rightarrow \infty} L_k$ is infinite, even
without assuming validity of DC. This can be seen using
an easy bound $R_r(k) > (k-1)^r$ implied by results
obtained by Abbott \cite{Abb65}
(cf. (7) and (4) in \cite{XXER}).

\bigskip
Observe an obvious equivalence that
$R(s,t) \geq R(s-1, t+1)$ if and only if
$R(s,t) - R(s-1, t) \geq R(s-1, t+1) - R(s-1, t)$,
which
for $3 \leq s \leq t$ can be seen as just another
way of looking at the DC.
It might seem that the analysis of how $R(s,t)-R(s-1,t)$
relates to $R(s-1,t+1)-R(s-1,t)$ should be simpler, but
it apparently resists to be so.
Some related discussion can be found
in \cite{BEFS,XSR,XXR}.

\bigskip
\section{Current Evidence for DC}

This section presents some additional observations
which make us believe that DC holds.
We note that Wang Rui \cite{Wang} did not provide
much intuition behind the conjecture itself, perhaps
because he thought that he had proved it as a theorem.
If so, then more discussion would not be required.

Below, we split our comments into two cases: of two colors
and of more colors. Only just a few exact values of
Ramsey numbers are known, hence not many absolute
instances confirming the DC can be pointed to.
On the other hand, for a large number of open cases, say
such as $R(s,t)$ for specific $s$ and $t$, it seems that
the best known lower bound is much closer to the exact
value than known upper bound. Historically (see the past
revisions of \cite{Radziszowski}), the lower bounds
often slowly improve over some time then stabilize,
while the upper bounds are improved rarely and most of
the time only with a large computational effort.
Or, in other words, known upper bounds are far from
being tight because we know very little
about how to improve them.
Thus, similarly as in the previous section when
arguing for the finiteness of 
$\lim_{r \rightarrow \infty} R_r(3)^{1 \over r}$,
our evidence will rely greatly on what we know about
lower bounds.

\medskip
\subsection*{Two Colors}

Let $DC(s,t)$ stand for the validity of $R(s,t) \geq R(s-1, t+1)$.
We will consider various $DC(s,t)$ statements for special values
of the parameters
$s$ and $t$, but always satisfying $3 \le s \le t$.

\begin{enumerate}[(a)]
\item
$DC(3,t)$ is true, since easily
$R(3,t) > R(2,t+1)=t+1$ for all $t \geq 3$.
\item
$DC(4,t)$ is true, since we have
$R(4,t) \geq R(3,t)+2t-3$ (\cite{BEFS}, see also \cite{XXR})
but easily $R(3,t+1) \leq R(3,t)+t+1$ for all $t \geq 4$.
\item
$DC(5,5)$ is true, since it is known that $R(5,5) \geq 43$
and $R(4,6) \le 41$ (cf. \cite{Radziszowski}).
Angeltveit and McKay in a recent unpublished project
\cite{AnMc} obtained the upper bounds $R(4,7) \le 58$ and
$R(4,8) \le 79$, which confirm the validity of
$DC(5,6)$ and $DC(5,7)$ by using previously
published lower bounds $R(5,6) \geq 58$ and
$R(5,7) \geq 80$ (cf. \cite{Radziszowski}).
\item
The above establishes the validity of $DC(s,t)$ for all
$s<5$ and all cases with $s+t \le 12$, except
$DC(6,6)$. For the latter it is known that
$R(6,6) \ge 102$ and $R(5,7) \le 143$, though recall
our previous comments, and especially in this case
we feel that the lower bound is strong but the upper
bound very weak.

In one larger case, namely that of $DC(8,10)$, 
the best known lower bounds
$R(8,10) \ge 343$ and $R(7,11) \ge 405$
(cf. \cite{Radziszowski}) do not "follow" the
DC (but do not contradict it either). 
We believe that this is because of a rather special
construction establishing the bound for $R(7,11)$, while
the bound for $R(8,10)$ was obtained by a
heuristic search restricted to only circular graphs.
This suggests that it should be feasible to
significantly improve the current lower bound
for $R(8,10)$.

We also note that
known bounds for $R(s,t)$ collected in \cite{Radziszowski}
do not contradict $DC(s,t)$ for any $3 \le s \le t$.

\item
The further we go from the diagonal of DC,
the easier it seems to corroborate it.
We anticipate this problem to be the hardest
on the diagonal itself, i.e. proving that
$R(t,t) \geq R(t-1,t+1)$ for any $t \ge 6$.

\item
In 2010,
Bohman and Keevash \cite{BK2010} proved that
for fixed $s \geq 5$ and $t \rightarrow \infty$ we have
the following lower bound
$$R(s,t) = \Omega(t^{{s+1} \over 2}
(\log s)^{{1 \over {s-2}}- {{s+1} \over 2}}).$$

This result does not resolve any concrete $DC(s,t)$
instances, yet, again using our perspective on lower
bounds, builds up evidence for the validity of
$DC(s,t)$ for fixed $s$ and large $t$.
\end{enumerate}

\medskip
\subsection*{More Colors}

In the multicolor cases, almost all evidence we have
for DC is based on lower bounds, even more so than in
the case of two colors. Table 2 lists 11 pairs
of parameters $(P_1,P_2)$ together
with the corresponding
best known lower bounds $(LB_1,LB_2)$ listed
in \cite{Radziszowski} for
$R(k_1, \cdots, k_{r-2}, k_{r-1}-1, k_r +1)$ and
$R(k_1, \cdots, k_r)$, with $4 \le k_{r-1} \le k_r$.
This includes essentially all
evidence of this type we have for $k_{r-1} \ge 4$.

\medskip
\begin{center}
\begin{tabular}{c|c||c|c}
$P_1$&$LB_1$&$LB_2$&$P_2$\cr
\hline
 3,3,5  &45 & 55 & 3,4,4\cr
 3,3,6  &61 & 89 & 3,4,5\cr
 3,3,7  &85 & 117 & 3,4,6\cr
 3,3,8  &103 & 152 & 3,4,7\cr
 3,3,9  &129 & 193 & 3,4,8\cr
 3,3,10  &150 & 242 & 3,4,9\cr
 3,4,6  &117 & 139 & 3,5,5\cr
 3,4,7  &152 & 181 & 3,5,6\cr
 3,4,8  &193 & 241 & 3,5,7\cr
 4,3,5  &89 & 128 & 4,4,4\cr
 3,3,3,5  &162 & 171 & 3,3,4,4\cr
\hline
\end{tabular}

\bigskip
{\bf Table 2.} The best known lower bounds $LB_1$ and $LB_2$
on Ramsey numbers
$R(k_1, \cdots, k_{r-2}, k_{r-1}-1, k_r +1)$
and $R(k_1, \cdots, k_r)$,
for some DC-adjacent pairs $(P_1,P_2)$, where
$P_1=(k_1, \cdots, k_{r-2}, k_{r-1}-1,k_r +1)$ and
$P_2=(k_1, \cdots, k_r)$.
\end{center}

\bigskip
We can say a little more beyond Table 2 for some
combinations of parameters in $P_2$ involving $k_{r-1}=3$.
For example, we clearly have $R(5,k)=R(k,2,5)$,
and by inspection of bounds reported
in \cite{Radziszowski}, we can see that
$R(k,3,4) \geq R(k,2,5)$ holds for $2 \le k \le 7$.

The lower bounds in columns $LB_1$ and $LB_2$ do get
occasional improvements, though not often and not by much.
For a particular $P_1$ to contradict DC, the corresponding
lower bound $LB_1$ would have to exceed not only $LB_2$,
but also its associated upper bound.

\bigskip
\section{Some Problems Related to DC and $R_r(k)$}

{\bf (1)} For connected graphs $G_1, \cdots, G_r$,
the generalized multicolor Ramsey number
$R(G_1, \cdots, G_r)$ is defined as the smallest
integer $n$ such that in any $r$-coloring
of the edges of $K_n$ there must be a monochromatic
$G_i$ in color $i$, for some $1 \le i \le r$.
We pose the following question
generalizing DC. For $G_{r-1}=K_s$,
$G_r=K_t$ with $s \le t$, is it true that
\vspace{-.1cm}
$$R(G_1, G_2, \cdots,
K_{s-1},K_{t+1}) \leq R(G_1,G_2, \cdots, K_s,K_t)?$$

\vspace{-.2cm}
\noindent
We think that it is true, but stop here
and do not make it another conjecture.

\bigskip
\noindent
{\bf (2)} Let $r \geq 3$, $k_i \ge 3$, and $k_{r-1} \leq k_r$.
Suppose that $C$ is a coloring of the edges of $K_n$
witnessing the lower bound
$n<R(k_1, \cdots, k_{r-2}, k_{r-1}-1, k_r +1)$.
Define the graph $G$ to consist of the edges of $C$
in colors $r-1$ and $r$. Is it true that
$G \not\rightarrow (k_{r-1}, k_r)^e$, i.e. that
there exists a 2-coloring of the edges of $G$
without any monochromatic $K_{k_{r-1}}$ in
the first color and $K_{k_r}$ in the second color?
We think that the answer is YES, but less strongly
than in (1).

\bigskip
\noindent
{\bf (3)} The Shannon capacity
of a noisy channel modeled by graph $G$, often
referred to as the Shannon capacity of $G$, is
defined as the limit
\vspace{-.1cm}
$$c(G) = \lim_{r \rightarrow \infty} \alpha(G^r)^{1 \over r},$$

\vspace{-.2cm}
\noindent
where $\alpha(G^r)$ is the independence number of
the strong $r$-th power of $G$.
The capacity $c(G)$ measures the efficiency of the best
possible strategy when sending long words over a noisy
channel modeled by $G$. It was studied extensively in
information theory by many authors, including
\cite{Alon2006,Alon1995}.
In a very short 1971 paper, Erd\H{o}s et al.
\cite{Erdos1971} proved that for each $k$ there
exists a graph $G$ with $\alpha(G) = k$
such that $\alpha(G^r)+1 = R_r(k+1)$. This provides an
implicit link between Shannon capacity and Ramsey numbers,
and in particular to the problem of finiteness of
the limit $\lim_{r \rightarrow \infty}R_r(k)^{1/r}$.
We explored it further in \cite{XRShannon}, where we proved
that $\lim _{n \rightarrow \infty}R_r(3)^{1/r}$ is
the supremum of the Shannon capacity of complements
of $K_3$-free graphs
but it cannot be achieved by any finite graph power.
In general, for any fixed integer $k \geq 3$,
we have that
$\lim _{r \rightarrow \infty}R_r(k)^{1/r}$
is equal to the supremum of the Shannon capacity $c(G)$
over all graphs $G$ with independence number $k-1$,
but this supremum cannot be achieved
by any finite graph power either.

\bigskip
\section*{Acknowledgement}

The work of the first and the third authors was partially supported by
the National Natural Science Foundation award (11361008).

\end{document}